# A bivariate version of the spread polynomials

Johann Cigler

*Dedicated to the memory of Hans-Christian Herbig*

**Abstract.**

We study a bivariate variant of Hans-Christian Herbig's version of Norman Wildberger's spread polynomials.

## 1. Introduction

The spread polynomials $S_n(x)$ introduced by Norman Wildberger satisfy $S_n\left(\sin^2\theta\right)=\sin^2(n\theta)$. The normalized version $Z_n(x)=4S_n\left(\frac{x}{4}\right)$ which satisfies $Z_n\left(4\sin^2\theta\right)=4\sin^2(n\theta)$ has been considered by Hans-Christian Herbig [5] and studied in [3] and [4]. It is closely related to the well-known Fibonacci and Lucas polynomials. In this note we study a bivariate variant $Z_n(x,s)$.

## 2. Some background material

In order to simplify the exposition, let me first recall some basic properties of the bivariate Fibonacci polynomials $F_n(x,s)$ and Lucas polynomials $L_n(x,s)$. They are well known, but readers who are unfamiliar with them can find everything proven in the introductory section of [2].

The Fibonacci polynomials $F_n(x,s)$ satisfy $F_n(x,s)=xF_{n-1}(x,s)+sF_{n-2}(x,s)$ with initial values $F_0(x,s)=0$ and $F_1(x,s)=1$.

The Lucas polynomials $L_n(x,s)$ satisfy the same recursion $L_n(x,s)=xL_{n-1}(x,s)+sL_{n-2}(x,s)$, but with initial values $L_0(x,s)=2$ and $L_1(x,s)=x$.

Their generating functions are

(1) $$\sum_{n\geq 0}F_n(x,s)z^n=\frac{z}{1-xz-sz^2}=\frac{z}{1-xz}\frac{1}{1-\frac{sz^2}{1-xz}}$$

and

(2) $$\sum_{n\geq 0}L_n(x,s)z^n=\frac{2-xz}{1-xz-sz^2}=\frac{2-xz}{1-xz}\frac{1}{1-\frac{sz^2}{1-xz}}.$$

Binet's formulae give

(3)
$$F_n(x,s)=\frac{\gamma^n-\bar{\gamma}^n}{\gamma-\bar{\gamma}},$$
$$L_n(x,s)=\gamma^n+\bar{\gamma}^n$$

where

(4)
$$\gamma=\gamma(x,s)=\frac{x+\sqrt{x^2+4s}}{2},\quad \bar{\gamma}=\bar{\gamma}(x,s)=\frac{x-\sqrt{x^2+4s}}{2}$$

are the roots of the characteristic polynomial $z^2-xz-s$ of the recurrence.

It is also well known that

(5)
$$F_n(x,s)=\sum_{k=0}^{\left\lfloor\frac{n-1}{2}\right\rfloor}\binom{n-1-k}{k}s^k x^{n-1-2k},$$
$$L_n(x,s)=F_{n+1}(x,s)+sF_{n-1}(x,s)=\sum_{k=0}^{\left\lfloor\frac{n}{2}\right\rfloor}\binom{n-k}{k}\frac{n}{n-k}s^k x^{n-2k}$$

and that (cf. [2],(2.24))

(6)
$$\sum_{j=0}^{n}(-s)^j\binom{2n}{j}L_{2n-2j}(x,s)=x^{2n}+(-s)^n\binom{2n}{n},$$
$$\sum_{j=0}^{n}(-s)^j\binom{2n+1}{j}L_{2n+1-2j}(x,s)=x^{2n+1}.$$

We shall also need Cassini's identity

(7)
$$F_n^2(x,s)-F_{n-1}(x,s)F_{n+1}(x,s)=(-s)^{n-1}.$$

**3. The polynomials $Z_n(x,s)$**

If we set $l_n(x)=L_n(x,-1)$ and $F_n(x)=F_n(x,1)$ then (cf. [3], [4])

(8)
$$Z_n(x)=2-l_n(2-x)=(-1)^{n-1}xF_n^2\left(\sqrt{x-4}\right).$$

The first terms are

$$\left(Z_n(x)\right)_{n\geq 0}=\left(0,x,4x-x^2,9x-6x^2+x^3,16x-20x^2+8x^3-x^4,25x-50x^2+35x^3-10x^4+x^5,\cdots\right).$$

Since $l_{mn}(x)=l_m\left(l_n(x)\right)$ and $l_2(x)=x^2-2$ imply $l_{2n}(x)=l_n\left(l_2(x)\right)=l_n\left(x^2-2\right)$ (cf. [3], (2.6)) we can also write

(9)
$$Z_n(x)=(-1)^{n-1}\left(l_{2n}\left(\sqrt{x}\right)-2(-1)^n\right).$$

Formula (9) suggests defining $Z_n(x,s)$ by

(10) $$Z_n(x,s)=L_{2n}(\sqrt{x},s)-2s^n.$$

The first terms are

(11) $$\left(Z_n(x,s)\right)_{n\geq 0}=\left(0,x,4sx+x^2,9s^2x+6sx^2+x^3,16s^3x+20s^2x^2+8sx^3+x^4,\cdots\right).$$

By (10) we get

(12) $$Z_n(x,s)=\alpha^{2n}+\bar{\alpha}^{2n}-2s^n$$

with

(13) $$\alpha=\gamma\left(\sqrt{x},s\right)=\frac{\sqrt{x}+\sqrt{x+4s}}{2},$$
$$\bar{\alpha}=\bar{\gamma}\left(\sqrt{x},s\right)=\frac{\sqrt{x}-\sqrt{x+4s}}{2}.$$

Let $\beta=\gamma\left(\sqrt{x+4s},-s\right)$ and $\bar{\beta}=\bar{\gamma}\left(\sqrt{x+4s},-s\right)$ be the roots of $z^2-\sqrt{x+4s}z+s=0.$

Then $\alpha=\beta,\ \bar{\alpha}=-\bar{\beta}$ and $\beta\bar{\beta}=s.$

Since $\alpha^{2n}+\bar{\alpha}^{2n}-2s^n=\beta^{2n}+\bar{\beta}^{2n}-2s^n=\left(\beta^n-\bar{\beta}^n\right)^2=x\left(\frac{\beta^n-\bar{\beta}^n}{\beta-\bar{\beta}}\right)^2=xF_n\left(\sqrt{x+4s},-s\right)$

we get

(14) $$Z_n(x,s)=xF_n\left(\sqrt{x+4s},-s\right)^2.$$

Let us also note that

(15) $$Z_{2n+1}(x,s)=L_{2n+1}^2(\sqrt{x},s),$$
$$Z_{2n}(x,s)=(x+4s)F_{2n}^2\left(\sqrt{x},s\right).$$

This follows from

$\alpha^{2(2n+1)}+\bar{\alpha}^{2(2n+1)}-2s^{2n+1}=\left(\alpha^{2n+1}+\bar{\alpha}^{2n+1}\right)^2$ and

$\alpha^{2(2n)}+\bar{\alpha}^{2(2n)}-2s^{2n}=\left(\alpha^{2n}-\bar{\alpha}^{2n}\right)^2=\left(\frac{\alpha^{2n}-\bar{\alpha}^{2n}}{\alpha-\bar{\alpha}}\right)^2(x+4s).$

Since $\alpha$ and $\bar{\alpha}$ are the roots of $z^2-\sqrt{x}z-s=0$ we see from $\left(\alpha^2-s\right)^2=x\alpha^2$ that $\alpha^2$ is a root of $z^2-(x+2s)z+s^2=0.$

By (12) $Z_n(x,s)$ is a linear sum of $n-$th powers of $\alpha^2,\bar{\alpha}^2,s,$ which are roots of

$(z-s)\left(z^2-\left(x+2s\right)z+s^2\right)=z^3-\left(x+3s\right)z^2+s\left(x+3s\right)z-s^3=0.$

Therefore, the sequence $Z_n(x,s)$ satisfies the recurrence

(16) $$Z_{n+3}(x,s)-(x+3s)Z_{n+2}(x,s)+s(x+3s)Z_{n+1}(x,s)-s^3Z_n(x,s)=0.$$

More precisely we get

(17) $$\sum_{n\geq 0} Z_{n+1}(x,s)z^n=\frac{x(1+sz)}{1-(x+3s)z+s(x+3s)z^2-s^3z^3}=\frac{x(1+sz)}{(1-sz)^3}\frac{1}{1-\frac{xz}{(1-sz)^2}}.$$

Cassini's identity gives another recurrence

(18) $$Z_{n-1}(x,s)Z_{n+1}(x,s)=\left(Z_n(x,s)-s^{n-1}x\right)^2$$

because

$$Z_{n+1}(x,s)Z_{n-1}(x,s)=x^2F_{n+1}^2\left(\sqrt{x+4s},-s\right)F_{n+1}^2\left(\sqrt{x+4s},-s\right)$$
$$=x^2\left(F_n^2\left(\sqrt{x+4s},-s\right)^2-s^{n-1}\right)^2=\left(xF_n^2\left(\sqrt{x+4s},-s\right)^2-s^{n-1}x\right)^2.$$

By (6) we get

(19) $$\sum_{j=0}^{n}(-s)^j\binom{2n}{j}Z_{n-j}(x,s)=x^n.$$

(5) implies

$$L_{2n}(x,s)=\sum_{k=0}^{n}\binom{2n-k}{k}\frac{2n}{2n-k}x^{2n-2k}s^k=\sum_{k=0}^{n}\binom{n+n-k}{2n-2k}\frac{2n}{n+n-k}x^{2n-2k}s^k$$
$$=\sum_{\ell=0}^{n}\binom{n+\ell}{2\ell}\frac{2n}{n+\ell}x^{2\ell}s^{n-\ell}=2s^n+\sum_{\ell=1}^{n}\binom{n+\ell-1}{2\ell-1}\frac{n+\ell}{2\ell}\frac{2n}{n+\ell}x^{2\ell}s^{n-\ell}$$
$$=2s^n+\sum_{\ell=1}^{n}\binom{n+\ell-1}{n-\ell}\frac{n}{\ell}x^{2\ell}s^{n-\ell}.$$

Thus,

(20) $$Z_n(x,s)=\sum_{k=1}^{n}c(n,k)s^{n-k}x^k$$

with

(21) $$c(n,k)=\frac{n}{k}\binom{n+k-1}{n-k}=\frac{n}{k}\binom{n+k-1}{2k-1}=\binom{n+k}{2k}+\binom{n+k-1}{2k}$$
$$=\frac{2}{(2k)!}n^2\left(n^2-1^2\right)\left(n^2-2^2\right)\cdots\left(n^2-(k-1)^2\right).$$

For example,

(22)
$$\left(c(n,k)\right)_{n,k=1}^{5}=\begin{pmatrix}1&0&0&0&0\\4&1&0&0&0\\9&6&1&0&0\\16&20&8&1&0\\25&50&35&10&1\end{pmatrix}.$$

Let us note that

(23)
$$Z_n(x)=(-1)^{n-1}Z_n(x,-1)=-\sum_{k=1}^{n}c(n,k)(-x)^k,$$
$$Z_n(x,s)=-s^nZ_n\left(-\frac{x}{s}\right).$$

**Remark**

OEIS A156308 shows that the polynomials $\sum_{k=1}^{n}c(n,k)x^k$ have also occurred in other contexts:

Peter Bala considered $R(n,x)=\frac{2}{x}\left(T_n\left(\frac{x+2}{2}\right)-1\right)$, where $T_n(x)$ denotes the Chebyshev polynomial of the first kind. $2T_n(x)=l_n(2x)$ gives
$xR(n,x)=l_n(x+2)-2=-Z_n(-x)=Z_n(x,1)$ by (8).

Alexander Burstein and Louis W. Shapiro defined them in [1],(5.1), via property (15).